 \documentclass[12pt]{amsart}
 \usepackage{amsmath,amssymb}
 \usepackage{pstricks}
 \usepackage{pst-node,pst-tree}
 \newtheorem{thm}{Theorem}[section]

 \theoremstyle{definition}
 
 \theoremstyle{remark}
 
 \numberwithin{equation}{section}
 \renewcommand{\a}{\mathbf{a}}
 \renewcommand{\b}{\mathbf{b}}
 \renewcommand{\c}{\mathbf{c}}
 \renewcommand{\d}{\mathbf{d}}
 \newcommand{\x}{\mathbf{x}}
 \newcommand{\R}{\mathbb{R}}
 
 \newcommand{\conv}{\mathrm{conv\;}}
 \newcommand{\set}[1]{\left\{#1\right\}}

\begin{document}
 \title{The Polytope of Dual Degree Partitions}%
 \author{Amitava Bhattacharya}
 \address{Bhattacharya: The University of Illinois at Chicago, United States}
 \email{amitava@math.uic.edu}
 \thanks{$^1$ This problem was raised in the Laplace Energy group of the
   Workshop \emph{Spectra of Families of Matrices described by
   Graphs, Digraphs, and Sign Patterns} held at the American
   Institute of Mathematics Research Conference Center on October
   23--27, 2006~\cite{Bru:06}.}
 \author{Shmuel Friedland}
 \address{Friedland: The University of Illinois at Chicago, United States}
 \email{friedlan@uic.edu}
 \author{Uri N. Peled}
 \address{Peled: (Corresponding Author) The University of Illinois at Chicago,
 MSCS Dept.\ (M/C 249), 851 S. Morgan Street, Chicago, IL 60607-7045,
 United States, Tel.\ 312-413-2156, Fax 312-996-1491} \email{uripeled@uic.edu}
 \date{November 21, 2006}%
 \begin{abstract}
  We determine the extreme points and facets of the convex hull of
  all dual degree partitions of simple graphs on $n$
  vertices.~$^1$
 \end{abstract}
 \maketitle
 \section{Introduction}
  We deal throughout with simple graphs $G$ (undirected, no loops,
  no multiple edges) on the vertex set $\set{1,\ldots,n}$. The
  \emph{degree} $d_i$ of vertex $i$ is the number of neighbors of
  $i$, and the \emph{degree sequence} of $G$ is $\d =
  (d_1,\ldots,d_n)$. We assume that the vertices have been relabeled
  so that $n-1 \geq d_1 \geq \cdots \geq d_n \geq 0$, and to stress
  this we call $\d$ a \emph{degree partition}. The \emph{dual degree
  partition} of $G$ is the sequence $\d^* = (d^*_1,\ldots,d^*_n)$,
  where $d^*_j = \set{i : d_i \geq j}$, so that $n \geq d^*_1 \geq
  \cdots \geq d^*_n = 0$. Both $\d$ and $\d^*$ can be conveniently
  pictured as a \emph{Ferrers diagram}, which is an $n \times n$
  matrix of zeros and ones, where the ones in each row are to the
  left of the zeros, the ones in each column are above the zeros,
  the row sums are the $d_i$ and the column sums are the $d^*_i$
  (such a matrix is sometimes called \emph{maximal}).
  Following Berge, it is also convenient to use the
  \emph{corrected Ferrers diagram}, which is obtained from the
  Ferrers diagram by moving every one on the main diagonal to the
  end of its row, replacing it with a zero. The row sums of the
  corrected Ferrers diagram are of course the $d_i$; its column sums
  $\overline{d^*_i}$ (which are not necessarily in non-increasing
  order) are called the \emph{corrected conjugate degrees}, and we
  use the notation $\overline{\d^*} =
  (\overline{d^*_1},\ldots,\overline{d^*_n})$.
  Figure~\ref{fig:ferrers} illustrates these definitions.
  \begin{figure}[ht]
    \label{fig:ferrers}
    \begin{center}
     \[
        \begin{array}{c|cccccc}
           & 6 & 6 & 3 & 1 & 0 & 0 \\
         \hline
         4 & 1 & 1 & 1 & 1 & 0 & 0 \\
         3 & 1 & 1 & 1 & 0 & 0 & 0 \\
         3 & 1 & 1 & 1 & 0 & 0 & 0 \\
         2 & 1 & 1 & 0 & 0 & 0 & 0 \\
         2 & 1 & 1 & 0 & 0 & 0 & 0 \\
         2 & 1 & 1 & 0 & 0 & 0 & 0 \\
       \end{array}
     \qquad
        \begin{array}{c|cccccc}
           & 5 & 5 & 2 & 3 & 1 & 0 \\
         \hline
         4 & 0 & 1 & 1 & 1 & 1 & 0 \\
         3 & 1 & 0 & 1 & 1 & 0 & 0 \\
         3 & 1 & 1 & 0 & 1 & 0 & 0 \\
         2 & 1 & 1 & 0 & 0 & 0 & 0 \\
         2 & 1 & 1 & 0 & 0 & 0 & 0 \\
         2 & 1 & 1 & 0 & 0 & 0 & 0 \\
       \end{array}
     \]
    \end{center}
    \caption{Ferrers diagram (left) and corrected Ferrers diagram
    (right) representing $\d = (4,3,3,2,2,2)$, $\d^* =
    (6,6,3,1,0,0)$, $\overline{\d^*} = (5,5,2,3,1,0)$.}
  \end{figure}

  Not every maximal matrix is the Ferrers diagram of the degree
  partition of a simple graph. If it is, we say that its row and
  columns sums are \emph{realizable}. There are many criteria for
  realizability proved in~\cite{Mah:95}. The one we use here
  involves the relation of majorization of sequences. We say that a
  sequence $\a=(a_1,\ldots,a_n)$ \emph{majorizes} a sequence
  $\b=(b_1,\ldots,b_n)$, and write $\a \succcurlyeq \b$, if the sum
  $A_k$ of the largest $k$ components of $\a$ and the sum $B_k$ of
  the largest $k$ components of $\b$  satisfy $A_k \geq B_k$ for
  $k=1,\ldots,n$ with equality for $k=n$. we call the difference
  $A_k-B_k$ the $k$th \emph{slack}. The following theorem is
  well-known; see a proof in~\cite{Mah:95}.

  \begin{thm}[Berge]\label{thm:Berge}
    Consider an integral sequence $\d = (d_1,\ldots,d_n)$ with $n-1
    \geq d_1 \geq \cdots \geq d_n \geq 0$. Then $\d$ and $\d^*$ are
    realizable if and only if $\sum_{i=1}^n d_i$ is even and
    $\overline{\d^*} \succcurlyeq \d$.
  \end{thm}

  The convex hull of all realizable degree partitions of length $n$
  was studied in~\cite{Bha:06}, and in particular that paper
  determines its extreme points and facets. In this paper we do the
  same for the convex hull of all realizable dual degree partitions
  of length $n$.

 \section{Results}\label{sec:results}
  Since $d^*_n=0$ for every realizable $\d^*$ of length $n$, we will
  suppress it and consider the convex hull of dual degree sequences
  as a subset of $\R^{n-1}$. We treat separately the case of $n$
  even and the somewhat harder case of $n$ odd. We use similar
  techniques in both cases.

  \subsection{$n$ even}\label{subsec:even}
   Consider the $n$ points
   \begin{equation}\label{ak}
    \a^{(k)} =
    (\underbrace{n,\ldots,n}_{k},\underbrace{0,\ldots,0}_{n-1-k}),
    \quad k=0,\ldots,n-1.
   \end{equation}
   \begin{thm}\label{thm:neven}
     For even $n$, the facet-defining inequalities of the convex
     hull of the realizable dual degree partitions on $n$ vertices
     are
     \begin{equation}\label{facesneven}
       n \geq x_1 \geq x_2 \geq \cdots \geq x_{n-1} \geq 0.
     \end{equation}
     Its extreme points are the $\a^{(k)}$ defined in~(\ref{ak}).
   \end{thm}
   This theorem was proved independently in the Laplace Energy group
   of the AIM workshop~\cite{Bru:06}, and by the first and third
   authors of this paper.
   \begin{proof}
    Let $P = \conv{\set{\a^{(0)},\ldots,\a^{(n-1)}}}$, let $Q$ be
    the convex hull of the realizable dual degree partitions on
    $n$ vertices, and let $R$ be the polytope defined
    by~(\ref{facesneven}). Each $\a^{(k)}$  is a realizable dual
    degree partition. This can be seen for example from the fact
    that by K\"{o}nig's theorem, $K_{\frac{n}{2},\frac{n}{2}}$ is
    the union of $\frac{n}{2}$ edge-disjoint matchings. The union
    of $k$ of these matchings is $k$-regular for $0 \leq k \leq
    \frac{n}{2}$, and the complements of these graphs are
    $k$-regular for $\frac{n}{2} \leq k \leq n-1$. It follows that
    $\a^{(k)} \in Q$ and thus $P \subseteq Q$. Obviously $Q
    \subseteq R$. We also have that every extreme point $\x$ of
    $R$ is one of the $\a^{(k)}$ and therefore $R \subseteq P$.
    Indeed, since $R \subseteq \R^{n-1}$ is defined by the $n$
    linear inequalities~(\ref{facesneven}), these inequalities
    holds with equality at $\x$ with at most one exception, and
    therefore $\x$ is one of the $\a^{(k)}$. This proves that
    $P=Q=R$. Furthermore, this polytope is full-dimensional
    because its $n$ extreme points are affinely independent. Since
    none of the inequalities~(\ref{facesneven}) is a consequence
    of the others, these inequalities are its unique
    facet-defining inequalities
   \end{proof}

  \subsection{$n$ odd}\label{subsec:odd}
   Consider the $\frac{n+1}{2}$ points
   \begin{equation}\label{akodd}
    \a^{(k)} =
    (\underbrace{n,\ldots,n}_{2k},\underbrace{0,\ldots,0}_{n-1-2k}),
    \quad 0 \leq 2k \leq n-1,
   \end{equation}
   the $\frac{n^2-1}{8}$ points
   \begin{equation}\label{bkl}
    \b^{(k,l)} =
    (\underbrace{n,\ldots,n}_{2k},\underbrace{n-1,\ldots,n-1}_{2l+1},\underbrace{0,\ldots,0}_{n-2-2k-2l}),
    \quad 0 \leq 2k +2l \leq n-3,
   \end{equation}
   and the $\frac{n^2-1}{8}$ points
   \begin{equation}\label{ckl}
    \c^{(k,l)} =
    (\underbrace{n,\ldots,n}_{2k+1},\underbrace{1,\ldots,1}_{2l+1},\underbrace{0,\ldots,0}_{n-3-2k-2l}),
    \quad 0 \leq 2k +2l \leq n-3,
   \end{equation}
   altogether $\frac{(n+1)^2}{4}$ points.

   In analogy with Theorem~\ref{thm:neven}, we have the following
   result.

   \begin{thm}\label{thm:nodd}
     For odd $n$, the facet-defining inequalities of the convex
     hull of the realizable dual degree partitions on $n$ vertices
     the $n$ inequalities~(\ref{facesneven}) as well as the inequality
     \begin{gather}
       \label{facesnodd2}
       (x_1-x_2) + (x_3-x_4) + \cdots + (x_{n-2} - x_{n-1}) \leq n-1.
     \end{gather}
     Its extreme points are the $\a^{(k)}$, $\b^{(k,l)}$ and
     $\c^{(k,l)}$ given in~(\ref{akodd}), (\ref{bkl})
     and (\ref{ckl}).
   \end{thm}

   \begin{proof}
    As before, let $P$ be the convex hull of the points $\a^{(k)}$,
    $\b^{(k,l)}$ and $\c^{(k,l)}$, let $Q$ be the convex hull of the
    realizable dual degree partitions on $n$ vertices, and let $R$
    be the polytope defined by~(\ref{facesneven}) and
    (\ref{facesnodd2}).

    We will show that each of the points $\a^{(k)}$, $\b^{(k,l)}$
    and $\c^{(k,l)}$ is a realizable dual degree partition, and
    consequently $P \subseteq Q$. We do this using
    Theorem~\ref{thm:Berge}. If $\d^* = \a^{(k)}$, then $\d =
    (\underbrace{2k,\ldots,2k}_n)$ and $\overline{\d^*} =
    (\underbrace{n-1,\ldots,n-1}_{2k},\underbrace{2k}_1)$. When we
    consider the majorization inequalities $\overline{\d^*}
    \succcurlyeq \d$, each of the first $2k$ inequalities adds
    $n-1-2k \geq 0$ to the slack, and the next inequality leaves the
    slack unchanged and exhausts $\d^*$.
    \\
    If $\d^* = \b^{(k,l)}$, then $\d =
    (\underbrace{2k+2l+1,\ldots,2k+2l+1}_{n-1},\underbrace{2k}_1)$
    and $\overline{\d^*} =
    (\underbrace{n-1,\ldots,n-1}_{2k},\underbrace{n-2,\ldots,n-2}_{2l+1},
    \underbrace{2k+2l+1}_1,\underbrace{0,\ldots,0}_{n-2k-2l-2})$.
    Each of the first $2k$ majorization inequalities adds $n-2k-2l-2
    \geq 1$ to the slack, each of the next $2l+1$ inequalities adds
    $n-2k-2l-3 \geq 0$ to the slack, and next inequality exhausts
    $\d^*$.
    \\
    If $\d^* = \c^{(k,l)}$, then $\d =
    (\underbrace{2k+2l+2}_{1},\underbrace{2k+1,\ldots,2k+1}_{n-1})$
    and $\overline{\d^*} =
    (\underbrace{n-1,\ldots,n-1}_{2k+1},\underbrace{2k+1}_1,\underbrace{1,\ldots,1}_{2l+1},
    \underbrace{0,\ldots,0}_{n-2k-2l-3})$. The first majorization
    inequality adds $n-2k-2l-3 \geq 0$ to the slack. Each of the
    next $2k$ inequalities adds $n-1-(2k+1) \geq 2l+1$ to the slack,
    for a total slack of $(n-2k-2l-3)+2k(n-1-(2k+1)) \geq
    2k(2l+1)$. Each of the next $2l+1$ inequalities subtracts $2k$
    from the slack, which keeps it nonnegative, and exhausts
    $\overline{\d^*}$.

    Each realizable dual degree partition $\x$ obviously
    satisfies~(\ref{facesneven}). It also satisfies~(\ref{facesnodd2})
    because $(x_1-x_2) + (x_3-x_4) + \cdots + (x_{n-2} - x_{n-1})$
    is the number of vertices of odd degree, which is even, and $n$
    is odd. Consequently $Q \subseteq R$.

    We will show that each extreme point $\x$ of $R$ is one of the
    points $\a^{(k)}$, $\b^{(k,l)}$ and $\c^{(k,l)}$, and
    consequently $R \subseteq P$. Since the polytope $R \subseteq
    \R^{n-1}$ is defined by $n+1$ inequalities, at least $n-1$ of
    these inequalities hold with equality at $\x$ and at most two
    are strict. Obviously at least one of the
    inequalities~(\ref{facesneven}) must be strict. Suppose exactly one
    inequality in~(\ref{facesneven}) is strict. Then $n = x_1 = \cdots = x_p
    > x_{p+1} = \cdots = x_{n-1} = 0$ for some $0 \leq p \leq n-1$.
    By~(\ref{facesnodd2}) $p$ must be even, say $p=2k$. Therefore
    $\x = \a^{(k)}$. We may thus assume that exactly two
    inequalities of~(\ref{facesneven}) are strict, and~(\ref{facesnodd2})
    holds with equality. Then $n = x_1 = \cdots = x_p > x_{p+1} =
    \cdots = x_q > x_{q+1} = \cdots = x_{n-1} = 0$ for some $0 \leq
    p < q \leq n-1$. In other words, $\x =
    (\underbrace{n,\ldots,n}_{p},\underbrace{x_{p+1},\ldots,x_{p+1}}_{q-p},
    \underbrace{0,\ldots,0}_{n-1-q})$. If $p$ is even, say $p=2k$,
    then since~(\ref{facesnodd2}) holds with equality, $q-p$ must be
    odd, say $q-p = 2l+1$, and $x_{p+1} = n-1$. Therefore $\x =
    \b^{(k,l)}$. If $p$ is odd, say $p=2k+1$, then by the same
    reason $q-p$ is odd, say $q-p = 2l+1$, and $x_{p+1} = 1$.
    Therefore $\x = \c^{(k,l)}$.

    We have shown that $P \subseteq Q \subseteq R \subseteq P$, so
    $P=Q=R$. Once again, this polytope is full-dimensional since
    it contains the $n$ affinely independent points $\a^{(k)}$ and
    $\b^{(0,l)}$. Since none of the inequalities~(\ref{facesneven}) and
    (\ref{facesnodd2}) is a consequence of the others, they are
    the unique facet-defining inequalities of that polytope.

   \end{proof}

   An integral point satisfying the
   inequalities~(\ref{facesneven}) and (\ref{facesnodd2}) need not be a
   realizable dual degree partition even if the sum of its
   components is even. An example for $n=7$ is given by
   $\d^*=(5,3,3,3,3,3)$, which satisfies~(\ref{facesneven}) and
   (\ref{facesnodd2}), yet $\d=(6,6,6,1,1,0,0)$ is not realizable.
   Therefore to characterize realizable dual degree partitions we
   need nonlinear constraints.

 \end{document}